\documentclass[journal,twoside]{IEEEtran}

\usepackage{algorithm, algorithmic}

\ifCLASSINFOpdf
\else
   \usepackage[dvips]{graphicx}
\fi
\usepackage{url}
\usepackage{graphicx}
\usepackage{mathtools}

\usepackage{array,amssymb,bm,paralist}
\usepackage[]{amsmath}
\usepackage[utf8]{inputenc}
\usepackage[english]{babel}
\usepackage{bbm}
\usepackage{float}
\usepackage[normalem]{ulem} %provides strike-through

\def\Vecq{\mathbf{q}}

\def\Vecu{\mathbf{u}}
\def\Vecv{\mathbf{v}}

\def\Vecx{\mathbf{x}}
\def\Vecy{\mathbf{y}}
\def\Vecz{\mathbf{z}}

\DeclareMathOperator{\Tr}{Tr}

\def\VecA{\mathbf{A}}
\def\VecB{\mathbf{B}}

\def\VecD{\mathbf{D}}

\def\VecI{\mathbf{I}}

\def\VecS{\mathbf{S}}

\def\VecX{\mathbf{X}}

\DeclareMathOperator*{\argmin}{arg\,min} 
\DeclareMathOperator{\diag}{diag} 
\DeclareMathOperator{\sign}{sign} 

\def\transpose{^{\!\mathsf{T}}}

\def\fa{\; \; \forall \;}

\newtheorem{prop}{Proposition}

\usepackage{xcolor}

\begin{document}

\title{SWAGGER: Sparsity Within and Across Groups for General Estimation and Recovery}

\author{Charles Saunders
%, \emph{Student Member, IEEE},
and Vivek K Goyal
%, \emph{Fellow, IEEE}
%\thanks{This work was supported in part by
%the Defense Advanced Research Projects Agency REVEAL Program under contract number HR0011-16-C-0030.}
%\thanks{C. Saunders and V. K. Goyal are with the Department of Electrical and Computer Engineering, Boston University, Boston, MA 02215 USA (cs13@bu.edu; v.goyal@ieee.org).}
}

\markboth{}{Saunders \& Goyal: SWAGGER}
\maketitle

\begin{abstract}
Penalty functions or regularization terms that promote structured solutions to optimization problems are of great interest in many fields. Proposed in this work is a nonconvex structured sparsity penalty that promotes one-sparsity within arbitrary overlapping groups in a vector. This allows one to enforce mutual exclusivity between components within solutions to optimization problems. We show multiple example use cases (including a total variation variant), demonstrate synergy between it and other regularizers, and propose an algorithm to efficiently solve problems regularized or constrained by the proposed penalty.
\end{abstract}


\IEEEpeerreviewmaketitle

\def\figurename{Fig.}

\section{Introduction}

\IEEEPARstart{S}{parsity}-promoting regularization is an integral component of modern-day signal processing, optimization and machine learning.
The most prevalent method,
LASSO~\cite{LASSO}, uses the $\ell_1$ norm to promote solutions to optimization problems that have few nonzero coefficients.
More recently, within some application domains there has been significant interest in structured sparsity.
For instance, group sparsity aims to set entire groups of components within a vector entirely to zero and allow the components in other groups to take on any value~\cite{GLASSOthesis}.
This is often achieved by regularizing with the $\ell_{2,1}$ norm over groups,
called group LASSO (\mbox{G-LASSO})~\cite{GLASSO}.
%\VKGcomment{Is the G-LASSO problem convex?  The mention of convexity with E-LASSO below draws attention as possibly highlighting a distinction with G-LASSO, but not explicitly enough one way or the other.}
Conversely, there is interest in sparsity penalties that promote intra-group sparsity, that is, ensuring groups within a vector are themselves sparse. Such penalties have been used to great success in areas such as deep learning~\cite{NNsparsity}, computer vision~\cite{computervision1,computervision2}, and medicine~\cite{medicineEEG,medicinefMRI}.
Previous work in this area includes elitist~\cite{ELASSO1,ELASSO2}, or exclusive LASSO (E-LASSO)~\cite{ELASSO3}
formulations that are convex and result in sparse, disjoint groups, and the nonconvex ‘sparsity within and across groups’ (SWAG)~\cite{SWAG1} plus extensions~\cite{SWAG2} that allow for overlapping groups of components.
%\VKGcomment{The edits above are an attempt to be careful about LASSO (and its analogues) being an optimization problem rather than a `penalty'.  But I chose to not be as pedantic about `sparsity' vs.\ `sparsity-promoting'.  You can change my mind on all of this.}

We propose here a more general, nonconvex structured sparsity penalty that allows for sparsity within and across overlapping groups for general estimation and recovery (SWAGGER).
The SWAGGER formulation encodes mutual exclusivity between pairs of components, or a transform of the components, using an easily constructed sparsity structure matrix.
This results in one-sparse groups with minimal bias in the nonzero entries, where bias is typically an unwanted byproduct of using convex sparsity penalties such as elitist LASSO\@.
We demonstrate the utility of SWAGGER in a number of settings, including a novel total variation-style denoising in one and two dimensions and modeling occlusions plus allowing mutually exclusive discretizations in imaging problems.
Additionally, we introduce a novel algorithm to efficiently solve problems regularized or constrained by this proposed penalty. 

\section{Proposed penalty}

We introduce a general structured sparsity penalty
\begin{equation}
 R(\Vecx) = \Phi(\VecB \Vecx) \transpose \VecS \Phi(\VecB \Vecx),
\label{eqn:IQ}
\end{equation}
for $\Vecx \in \mathbb{R}^N$.
For certain choices of
the matrix $\VecB \in \mathbb{R}^{M \times N}$, where $M$ is any integer,
function $\Phi : \mathbb{R}^M \rightarrow \mathbb{R}^M$,
and matrix $\VecS \in \mathbb{R}^{M \times M}$,
this penalizes deviation from a selected form of structured sparsity.
Possible prototypes for $\VecB$ include the identity matrix, finite difference matrices, or a transform to a different basis; see Section~\ref{sec:examples}.
The function $\Phi$ introduces a nonlinearity that helps to bound the penalty from below. We constrain $\Phi$ in Prop.~\ref{prop:constraints} below, and examples are discussed in Section~\ref{sec:phichoice}.
The matrix $\VecS$ is key, as it encodes the structural sparsity properties. We define $\VecS$ as:
\begin{equation}
    \VecS_{i,j} = \VecS_{j,i} = \mu,
    \label{eqn:Sdef}
\end{equation}
%\VecS_{i,j} = \VecS_{j,i} = P(\Vecz_j = 0 \, | \, \Vecz_i \neq 0)$,
where $\mu \in [0,1]$ is the \emph{strength of exclusivity} between a pair of components, $\Phi(\VecB \Vecx)_i$ and $\Phi(\VecB \Vecx)_j$.
When $\mu = 0$, there is no exclusivity between the two components.
When $\mu > 0$, it indicates that the two components should not be nonzero simultaneously (with the strength of exclusivity increasing with $\mu$).
We bound $\mu$ in the [0,1] range for convenience and note that any scaling can be absorbed into a multiplicative regularization parameter.
Given this interpretation, the diagonal of $\VecS$ will be zero and hence $\VecS$ will typically be indefinite, resulting in nonconvexity. Examples for $\VecS$ are given in Section~\ref{sec:Schoice}.

In order for the penalty function to be useful in regularizing minimization problems, it must be bounded from below.
Prop.~\ref{prop:constraints} provides simple sufficient conditions for this (proof in Appendix~\ref{sec:penbounds}).

\medskip

\begin{prop}
\label{prop:constraints}
In order for \eqref{eqn:IQ} to be bounded from below, it is sufficient to require:
\begin{enumerate}
    \item $\VecS_{i,j} \geq 0 \quad \fa (i,j) \in \{1,2,\ldots,M\}^2$; and
    \item $\Phi(\Vecx) \geq 0 \; \mbox{elementwise} \; \fa \Vecx$.
\end{enumerate}
\end{prop}

%\VKGcomment{I realize that by putting ``sufficient'' in the proposition, I am having this word repeated.
%Feel free to rephrase, but I found it troubling that the proposition seemed to be claiming necessity instead of sufficiency, and one shouldn't have to read the surrounding text to reinterpret.}

The definition for $\VecS$ in \eqref{eqn:Sdef} meets condition 1) in Prop.~\ref{prop:constraints}.
Except for a discussion of extensions in Section~\ref{sec:phichoice},
we assume $\Phi(\Vecx)$ is the absolute value $|\Vecx|$, which satisfies condition 2).

\section{Properties of the SWAGGER penalty}

%\VKGcomment{What do you think of swapping Sections III and IV?
%Or to move {\sc Algorithms} later, after {\sc Choices for $\Phi(\cdot)$}.
%That feels more natural to me, and I don't think it deemphasizes algorithms.}
%\CS{I didn't move the sections although I am open to it. I chose the current order because I view the algorithm as one of the main contributions of this paper so it might be good to have it up front, but I see why the other order makes sense..}

When $\Phi(\cdot) = |\cdot|$, the Hessian of the penalty is given by
\begin{equation*}
    \nabla^2_\Vecx R(\Vecx) = 2 \, \VecX \, \VecS \, \VecX,
\end{equation*}
where $\VecX = \diag(\sign(\Vecx))$.
The sparsity structure matrix, $\VecS$, is symmetric and has zeros along its main diagonal.
$\VecS$ must be indefinite as it has at least one positive eigenvalue (see~\cite{eigbound}) and
the sum of the eigenvalues is zero (since the trace of a symmetric matrix
is equal to the sum of its eigenvalues
and $\Tr(\VecS) = 0$).
The Hessian has the same structure and can therefore not be positive semidefinite outside of the trivial case where $\Vecx = 0$.
Thus, the penalty is \emph{always} nonconvex for practical purposes. 
More analysis of the eigenvalues of matrices of this kind can be found in~\cite{hollowmatrix}.

We can compare SWAGGER with the convex E-LASSO penalty, which uses the $\ell_1$-squared norm to achieve intra-group sparsity.
In the case of a one-sparse group (see Section~\ref{sec:onesparse}), we can introduce a modified sparsity matrix, $\VecS_\eta = \mathbbm{1}\mathbbm{1}\transpose - \eta \VecI$, where scalar parameter $\eta \in [0,1]$ gives us E-LASSO for $\eta = 0$ and SWAGGER for $\eta = 1 $:
\begin{equation*}
|\Vecx| \transpose \VecS_\eta  |\Vecx| = \|\Vecx\|_1^2 - \eta  \|\Vecx\|_2^2.
\end{equation*}
When $\eta =1$, we see the introduction of the $- \eta  \|\Vecx\|_2^2$ term which both de-biases the penalty and also ensures the penalty equals zero when $\Vecx$ is one-sparse. For a more general $\VecS$, the closest convex penalty would be to replace $\VecS$ with $\VecS - c \VecI$, where $c$ takes on the value of the smallest eigenvalue of $\VecS$.

\paragraph*{Maintaining convexity of a full cost function}
There has recently been interest in so called `convex nonconvex' sparsity regularization where parameters of a nonconvex penalty are set to values that maintain convexity of the full cost function with the data fidelity, e.g.,~\cite{Selesnick2017}.
This is achievable with the SWAGGER penalty, also.
Consider the problem of minimizing a cost function of the form
$\frac{1}{2} \| \VecA \Vecx - \Vecy \|_2^2 + \lambda |\Vecx| \transpose \VecS |\Vecx|$.
One can replace $\VecS$ with $\VecS - c\VecI$, where
$c = \lambda_\text{min}(\VecS) + \frac{1}{2 \lambda}  \lambda_\text{min}(\VecA \transpose \VecA)$,
to achieve nonconvex regularization whilst maintaining the overall convexity of the problem (see derivation in Appendix~\ref{sec:cnc}).
This is only useful when the smallest eigenvalue of $\VecA \transpose \VecA$ is not too small and the regularization strength parameter $\lambda$ is not too large.
The work presented outside of this section is only focused on the use of the unmodified $\VecS$, but further analysis of the intermediate convexity paradigm may be of interest.

\section{Notable cases for $\VecS$}
\label{sec:Schoice}

\subsection{Canonical One-sparsity}
\label{sec:onesparse}
\begin{figure}
\centerline{\includegraphics[width=\columnwidth]{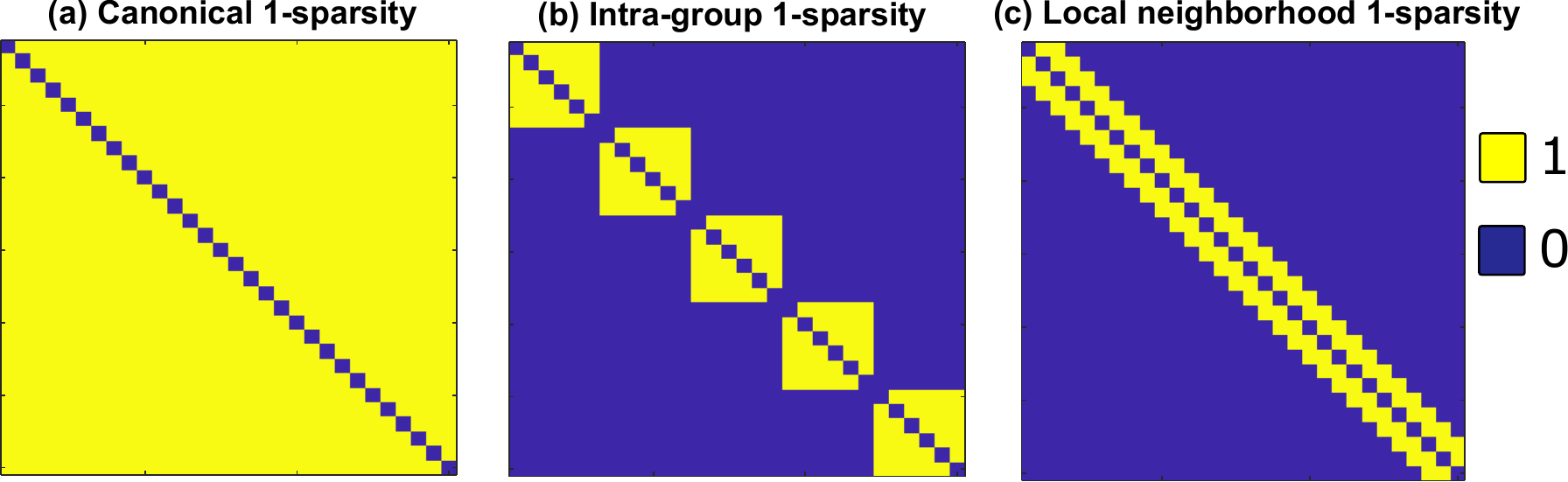}}
\caption{\textbf{(a)} $\VecS$ matrix that promotes a one-sparse solution. \textbf{(b)} $\VecS$ matrix that promotes one-sparse groups in the solution. \textbf{(c)} $\VecS$ matrix that promotes a minimum separation distance between nonzero components in the solution. }
\label{fig:Smatrix}
\end{figure}

Define $\Phi(\cdot) = |\cdot|$ and $\VecB = \VecI$.
In this case, we obtain
\begin{equation}
 R(\Vecx) = |\Vecx| \transpose \VecS |\Vecx|.
\label{eqn:IQAbs}
\end{equation}
This is the `sparsity within and across overlapping groups' formulation which is introduced in~\cite{SWAG2}.
If we chose
$\VecS = \mathbbm{1}\mathbbm{1}\transpose - \VecI$,
where $\mathbbm{1}$ is a column vector of 1s of appropriate dimension,
then we see that
\begin{equation}
\begin{split}
 R(\Vecx) & = |\Vecx| \transpose \VecS |\Vecx| = |\Vecx| \transpose (\mathbbm{1}\mathbbm{1}\transpose - \VecI)|\Vecx|  \\ & = (|\Vecx| \transpose \mathbbm{1})^2 - |\Vecx| \transpose |\Vecx| = \|\Vecx\|_1^2 - \|\Vecx\|_2^2,
 \end{split}
\label{eqn:IQSc}
\end{equation}
which is equal to zero for any one-sparse $\Vecx$. Hence, this can be used to promote one-sparse solutions to optimization problems. See Fig.~\ref{fig:Smatrix}(a) for a graphical depiction.

\subsection{Intra-group Sparsity}\label{sec:grpsparse}
Define $\Phi(\cdot) = |\cdot|$ and $\VecB = \VecI$. Consider a vector $\Vecx \in \mathbb{R}^{gn}$ where $g$ is a number of groups and $n$ is the number of components in each group.
Then, we can promote one-sparsity within each group using $\VecS \in \mathbb{R}^{gn \times gn}$ defined in blocks:
\begin{equation*}
\VecS = \begin{bmatrix} \widetilde{\VecS} & 0 & \cdots & 0 \\
                           0 & \widetilde{\VecS} & \cdots & 0 \\
                           \vdots & \vdots & \ddots & 0 \\
                           0 & 0 & 0 & \widetilde{\VecS}
          \end{bmatrix},
\end{equation*}
where $\widetilde{\VecS} = \mathbbm{1}\mathbbm{1}\transpose - \VecI_{n}$ (see Fig.~\ref{fig:Smatrix}(b)).
This is now a block-separable problem, equivalent to applying the penalty defined in Section~\ref{sec:onesparse} to different sections of $\Vecx$, i.e.,
$\sum_{k=1}^g |\Vecx_k|\transpose \widetilde{\VecS} |\Vecx_k|$,
where $k$ indexes blocks of $n$ components.
One use case for this formulation is in modeling occlusions in imaging problems, which we address with an example in Section~\ref{sec:erti}. 

\subsection{Local Neighborhood Sparsity}\label{sec:LN}
Sparsity can be promoted in overlapping local neighborhoods in an ordered sequence of indexes
by restricting the nonzeros to a band around the diagonal.
This gives a symmetric Toeplitz $\VecS$ (see Fig.~\ref{fig:Smatrix}(c)).
For instance, an $\VecS$ with neighbor distance $n=1$ is given by

\begin{equation*}
\VecS = \begin{bmatrix} 0 & 1 & 0 & 0 & \hdots \\  1 &  0 & 1 & 0 & \ddots \\ 0 &  1 & 0 & 1 & \ddots \\ 0 & 0 & 1 & 0 & \ddots  \\ \vdots & \ddots & \ddots & \ddots & \ddots    \end{bmatrix}.
\end{equation*}

For a general neighbor distance $n$, $\VecS$ is given by
\begin{equation*}
\VecS_{i,j} = \begin{cases} 
          \mu,  & |i - j| \leq n, \; i\neq j; \\
          0, & \text{otherwise},
       \end{cases}
\end{equation*}

where again $\mu \in [0,1]$ defines a strength of exclusivity. This imposes local-$n$-neighborhood one-sparsity, which ensures nonzeros in a solution occur with a minimum separation distance of $n$ indices.

\section{Choices for $\Phi(\cdot)$}
\label{sec:phichoice}
The typical choice for $\Phi(\cdot)$ is the absolute value $|\cdot|$.
This is the main focus of this paper and admits fast algorithms to result in one-sparse groups with a prescribed structure, with minimal bias.
However, we note that other choices can enjoy interesting properties.
For instance, using
$\Phi(\Vecx) = \text{max}(0,\Vecx)$
or
$\Phi(\Vecx) = \text{max}(0,-\Vecx)$
acts on positive or negative values only, respectively, and meets the condition outlined in Prop.~\ref{prop:constraints}.
Thus, combining two constraint terms,
\begin{equation}
 \begin{aligned}
\hat{\Vecx} =
& \argmin_{\Vecx} 
& & f(\Vecx) \\
& \text{s.t.  }
& & \text{max}(0,\Vecx) \transpose \, \VecS \, \text{max}(0,\Vecx) = 0 \\
& & & \text{max}(0,-\Vecx) \transpose \, \VecS \, \text{max}(0,-\Vecx) = 0
\end{aligned}
\label{eqn:poscorr}
\end{equation}
can result in a solution where only \emph{positive} correlation is penalized, i.e., a group may contain one positive and negative component, rather than being one-sparse.
In~\cite{NNsparsity}, exclusive sparsity regularization is used to reduce redundancy in the learned kernels in convolution neural networks.
The $\Phi(\cdot) = \text{max}(0,\cdot)$ formulation could function as a relaxed approach to this, allowing more degrees of freedom in the filter kernels. 

In a case where there can be a large disparity between the size of components within groups in the solution vector, one may wish to impose a nonlinear scaling and use $\Phi(\Vecx) = |\Vecx|^\kappa$ where $\kappa < 1$.
This reduces the effect of large but possibly incorrect components on other, potentially correct components. This can be seen in use in the two-dimensional local neighborhood total variation example in Section~\ref{sec:lntv}.

\section{Algorithms}

First, we note that when $\Phi(\cdot) = |\cdot|$ and $\VecB = \VecI$, the penalty can be expressed
\begin{equation}
 R(\Vecx) = |\Vecx| \transpose \VecS |\Vecx| = \|\Vecx\|_1^2 - \|\Vecx\|_2^2 - |\Vecx|\transpose(\mathbbm{1}\mathbbm{1}\transpose - \VecS - \VecI)|\Vecx|.
\label{eqn:IQ_2}
\end{equation}
An efficient algorithm to evaluate the proximal operator of the $\ell_1^2$ norm is available (see Appendix~\ref{sec:l12})~\cite{l12proof}. Hence, we can aim to efficiently solve any problem of the form:
\begin{equation}
 \begin{aligned}
\hat{\Vecx} =
& \argmin_{\Vecx} 
& & f(\Vecx) \\
& \text{s.t.  }
& & |\Vecx| \transpose \VecS |\Vecx| = 0
\end{aligned}
\label{eqn:optprob}
\end{equation}
using a proximal subgradient method (Algorithm~\ref{alg:swag}).
This simply requires us to evaluate the (sub)gradient of
$f(\Vecx) - \|\Vecx\|_2^2 - |\Vecx|\transpose(\mathbbm{1}\mathbbm{1}\transpose - \VecS - \VecI)|\Vecx|$
and the proximal operator of $\|\Vecx\|_1^2$.
In practice, we find the convergence speed of this algorithm is not negatively affected by the use of a subgradient of
$-|\Vecx|\transpose(\mathbbm{1}\mathbbm{1}\transpose - \VecS - \VecI)|\Vecx|$
as this term will always aim to push the components of $\Vecx$ \emph{away} from the discontinuity in the gradient at $\Vecx = 0$.
Furthermore, in cases with no overlapping groups, this term equals zero. 

As this penalty is nonconvex, it is prudent to start with an initial estimate of $\Vecx$ that minimizes the data fidelity.
Given this, Algorithm~\ref{alg:swag} outlines the proximal subgradient scheme used to approximately solve the problem in \eqref{eqn:optprob} by introducing a Lagrange multiplier $\lambda$.

 \begin{algorithm}
 \begin{algorithmic}[1]
 \renewcommand{\algorithmicrequire}{\textbf{Input:}}
 \renewcommand{\algorithmicensure}{\textbf{Output:}}
 \renewcommand{\algorithmicfor}{\textbf{while}}
 \REQUIRE Initial estimate $\Vecx^0 = \argmin_{\Vecx} f(\Vecx)$, $\lambda^0 = 0$, \\ $\overline{\VecS} = \mathbbm{1}\mathbbm{1}\transpose - \VecS $.
 \ENSURE  $\hat{\Vecx}$
% \\ \textit{Initialisation} :
%  \STATE first statement
% \\ \textit{LOOP Process}
  \FOR {not converged}
  \STATE $\Vecx^{k+1} =P_{\alpha \lambda^k}( \Vecx^{k} - \alpha(\nabla_x (f(\Vecx^k) - \lambda^k|\Vecx^k|\transpose \, \overline{\VecS} \, |\Vecx^k|))$
  \STATE $\lambda^{k+1} = \lambda^k + \alpha(|\Vecx^{k}| \transpose \VecS |\Vecx^{k}|)$
  \ENDFOR
 \RETURN $\hat{\Vecx} = \Vecx^{k+1}$ 
 \end{algorithmic} 
  \caption{Proximal subgradient method for SWAGGER constrained problems}
  \label{alg:swag}
 \end{algorithm}

In Algorithm~\ref{alg:swag},
$P_{\alpha\lambda^k}$ is the proximal operator associated with the $\ell_1^2$ term, $\|\cdot\|_1^2$, and $\alpha$ is a step size that can be fixed or computed dynamically.
There exists an accelerated proximal gradient algorithm that is guaranteed to converge to a stationary point even for nonconvex problems~\cite{NCAPG}, which could be used here also.
The algorithm for this is outlined in Appendix~\ref{apdx:ncapg}. 

When $\VecB$ is not the identity, we can instead solve a problem of the form:
\begin{equation}
 \begin{aligned}
\hat{\Vecx} =
& \argmin_{\Vecx} 
& & f(\Vecx)  \\
& \text{s.t.} 
& &  \VecB \Vecx = \Vecz \\ 
& & & |\Vecz| \transpose \VecS |\Vecz| = 0
\end{aligned}
\label{eqn:admm}
\end{equation}
using the alternating direction method of multipliers (ADMM) \cite{admm}, which also uses Algorithm~\ref{alg:swag} or an accelerated variation thereof at each iteration.

\section{Example applications}
\label{sec:examples}

\subsection{Comparison with Other Methods}
\label{sec:compare}

Three experiments were performed, where measurements $\Vecy$ were generated using $\Vecy = \VecA \Vecx$, with white Gaussian noise added to achieved a signal-to-noise radio (SNR) of 25 dB\@.
Matrix $\VecA \in \mathbb{R}^{25\times60}$ has random Gaussian entries realized for each trial.
The vector $\Vecx$ has entries $\pm U(0.5,1.5)$ and is made to fit a specific sparsity structure such that $|\Vecx| \transpose \VecS |\Vecx| \approx 0$ for some $\VecS$.
Details of how this is achieved are presented in Appendix~\ref{sec:genx}.
In the first experiment, $\VecS$ is fixed and enforces intra-group one-sparsity in ten groups (see Section~\ref{sec:grpsparse}), where
$\widetilde{\VecS} \in \mathbb{R}^{6 \times 6}$.
In the second experiment, a local neighborhood $\VecS$ (see Section~\ref{sec:LN}) is used with $n=4$.
In the third experiment, a new binary $\VecS$ is generated per trial with entries equal to zero or one with $50\%$ probability.
We compare the performance of the SWAGGER penalty with other typical sparsity regularizers, both structured and unstructured.
The problems solved are of the form
$\hat{\Vecx} = \argmin_{\Vecx} \frac{1}{2}\|\VecA\Vecx -y\|_2^2 + \lambda R(\Vecx)$,
where $R(\Vecx)$ is:
\begin{enumerate}
    \item  LASSO: $\|\Vecx\|_1$ 
    \item  p-shrinkage~\cite{pshrinkage}, with $p=0.5$ 
    \item  E-LASSO/$\ell_{1,2}$: $\sum_i\|\Vecx_{gi}\|_1^2$
\end{enumerate}
The E-LASSO penalty uses a sum of the $\ell_1$-squared term over groups.
We create one group per row in $\VecS$, where each group contains the nonzero components in each row (and the diagonal), i.e.,
$g_i = \{\, j \mid \, \VecS_{i,j} + \delta_{i-j} = 1 \,\}$.
This penalty is not designed for overlapping groups, so the results are poor in the local neighborhood and random $\VecS$ experiments.
The pseudo-inverse is used in the initialization for the nonconvex penalties, $\Vecx^0 = (\VecA\transpose\VecA)^{-1}\VecA\transpose \Vecy$.

\newcommand{\indentA}{\hspace*{2mm}}
\newcommand{\indentB}{\hspace*{6mm}}

\begin{table}
\caption{Simulated results over 1000 trials for three different sparsity structures.
Metrics are the percentage of nonzero components identified, the Jacard index (0 - 1), and the mean-squared error (MSE) where the ground truth is nonzero.
}
\label{table:testresults}
\begin{center}
\renewcommand{\arraystretch}{1.1}
\begin{tabular}{@{}lrrr@{}}
\hline
\textbf{Experiment}                                                       & \textbf{\begin{tabular}[c]{@{}r@{}}Support correct \\ (\%)\end{tabular}} & \textbf{Jacard Index}       & \textbf{\begin{tabular}[c]{@{}l@{}}MSE in \\ support\end{tabular}} \\ \hline
\textbf{Group S}                                                          & \textbf{}                                                                &                             &                                                                    \\ 
                                                                           \multicolumn{3}{l}{\indentA \textbf{$\lambda$ tuned for sparsity level}}                                                 &                                                                    \\ 
\indentB SWAGGER                                                                   & \textbf{75.69}                                                           & \textbf{0.646}              & \textbf{1.987}                                                     \\
\indentB E-LASSO                                                                   & 59.84                                                                    & 0.398                       & 4.634                                                              \\
\indentB $p$-shrinkage                                                             & 68.89                                                                    & 0.559                       & 2.564                                                              \\
\indentB LASSO                                                                     & 56.71                                                                    & 0.404                       & 3.813                                                              \\ 
                                                                          \multicolumn{3}{l}{\indentA \textbf{$\lambda$ tuned for Jacard index}}                                                 &                                                                    \\ 
\indentB SWAGGER                                                                   & 84.10                                                                    & \textbf{0.699}              & \textbf{1.633}                                                     \\
\indentB E-LASSO                                                                   & \textbf{87.38}                                                           & 0.509                       & 2.849                                                              \\
\indentB $p$-shrinkage                                                             & 77.82                                                                    & 0.638                       & 2.122                                                              \\
\indentB LASSO                                                                     & 73.91                                                                    & 0.492                       & 3.221                                                              \\ \hline
\multicolumn{2}{l}{\textbf{Local neighborhood S}}
%\begin{tabular}[c]{@{}l@{}}Local neighbor-\\ hood S\end{tabular} }
& \textbf{}                   &                                                                    \\ 
                                                                          \multicolumn{3}{l}{\indentA \textbf{$\lambda$ tuned for sparsity level}}                                                &                                                                    \\ 
\indentB SWAGGER                                                                   & \textbf{79.89}                                                           & \textbf{0.701}              & \textbf{1.653}                                                     \\
\indentB E-LASSO                                                                   & 60.37                                                                    & 0.241                       & 5.432                                                              \\
\indentB $p$-shrinkage                                                             & 68.02                                                                    & 0.545                       & 2.671                                                              \\
\indentB LASSO                                                                     & 56.46                                                                    & 0.402                       & 3.813                                                              \\ 
                                                                          \multicolumn{3}{l}{\indentA \textbf{$\lambda$ tuned for Jacard index}}                                                 &                                                                    \\ 
\indentB SWAGGER                                                                   & \textbf{85.18}                                                           & \textbf{0.754}              & \textbf{1.328}                                                     \\
\indentB E-LASSO                                                                   & 82.18                                                                    & 0.306                       & 5.278                                                              \\
\indentB $p$-shrinkage                                                             & 76.92                                                                    & 0.625                       & 2.187                                                              \\
\indentB LASSO                                                                     & 72.10                                                                    & 0.479                       & 3.318                                                              \\ \hline
\textbf{Random S}                                                         & \textbf{}                                                                & \textbf{}                   &                                                                    \\ 
                                                                          \multicolumn{3}{l}{\indentA \textbf{$\lambda$ tuned for sparsity level}}                                                 &                                                                    \\ 
\indentB SWAGGER                                                                   & \textbf{95.32}                                                           & \textbf{0.927}              & \textbf{0.171}                                                     \\
\indentB E-LASSO                                                                   & 56.20                                                                    & 0.147                       & 3.151                                                              \\
\indentB $p$-shrinkage                                                             & 89.68                                                                    & 0.845                       & 0.461                                                              \\
\indentB LASSO                                                                     & 70.90                                                                    & 0.563                       & 1.708                                                              \\ 
                                                                          \multicolumn{3}{l}{\indentA \textbf{$\lambda$ tuned for Jacard index}}                                                 &                                                                    \\ 
\indentB SWAGGER                                                                   & \textbf{95.84}                                                           & \textbf{0.932}              & \textbf{0.151}                                                     \\
\indentB E-LASSO                                                                   & 88.64                                                                    & 0.206                       & 3.173                                                              \\
\indentB $p$-shrinkage                                                             & 93.15                                                                    & 0.888                       & 0.327                                                              \\
\indentB LASSO                                                                     & 88.26                                                                    & 0.661                       & 1.173                                                             \\ \hline
\end{tabular}
\end{center}
\end{table}

Table~\ref{table:testresults} shows the results from the three experiments.
The metrics evaluated are the percentage of components correctly identified as nonzero, the Jacard index (the intersection of the true support and estimated support divided by the union), and the mean-squared error (MSE) for the components that are nonzero in the ground truth.
We see that SWAGGER outperforms the other algorithms in almost every case.
It is unsurprising that the unstructured sparsity regularizers (LASSO and p-shrinkage) perform especially poorly in the group $\VecS$ and local neighborhood settings as the true support is not especially sparse overall,
whereas SWAGGER is agnostic to the overall sparsity level.
Furthermore, the E-LASSO ($\ell_1^2$) regularizer performs poorly when groups are overlapping, whereas SWAGGER handles this gracefully and without any additional effort required.
We also see that the MSE within the correct support is particularly low with SWAGGER, as it implicitly de-biases the solutions (this can be readily seen from the $-\|\Vecx\|_2^2$ term in (\ref{eqn:IQSc})).
Of particular interest are the results where $\lambda$ is tuned to result in the correct sparsity level.
The SWAGGER results are particularly favorable in this setting, and this is the most typical use case; it is roughly equivalent to solving the constrained SWAGGER problem in (\ref{eqn:optprob}) and hence requires no user tuning of the regularization strength.

\subsection{Local Neighborhood Total Variation}
\label{sec:lntv}
\begin{figure}
\centerline{\includegraphics[width=\columnwidth]{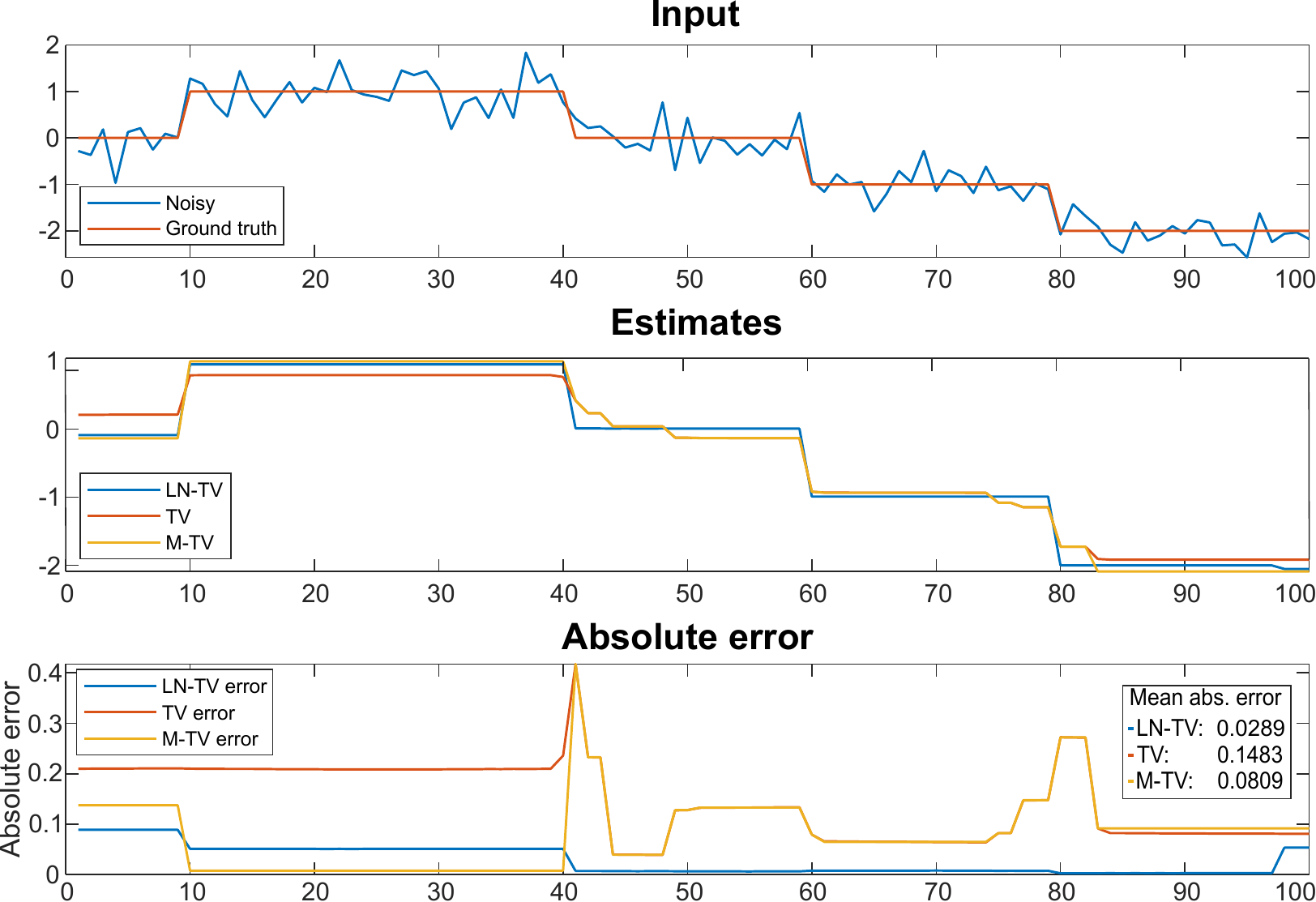}}
\caption{The proposed local neighborhood total variation (LN-TV) denoising, which only allows changes in the derivative to occur with some minimum separation, compared to standard total variation (TV) and Moreau-enhanced total variation (M-TV)  compared to  denoising.}
\label{fig:LNTVcompare}
\end{figure}

The well established total variation (TV) penalty, $\|\VecD \Vecx\|_1$ where multiplying by $\VecD$ takes finite differences, aims to sparsify the changes along the solution vector $\Vecx$.
We propose here a local neighborhood total variation (LN-TV) penalty using SWAGGER\@.
LN-TV is achieved by assigning $\VecB = \VecD$ and using an $\VecS$ that promotes sparsity in local neighborhoods as defined in Section~\ref{sec:LN}.
LN-TV simply ensures that changes in the vector occur at least some minimum distance apart, while allowing changes of any amplitude.
This is a prior which suggests the vector is piecewise-constant with a certain minimum segment length.
If a typical minimum separation distance between significant changes in a vector is known, LN-TV can outperform TV and Moreau-enhanced TV (a convex-nonconvex variation)~\cite{MTV} in a denoising setting, as seen in Fig.~\ref{fig:LNTVcompare}.

\begin{figure}
\centerline{\includegraphics[width=\columnwidth]{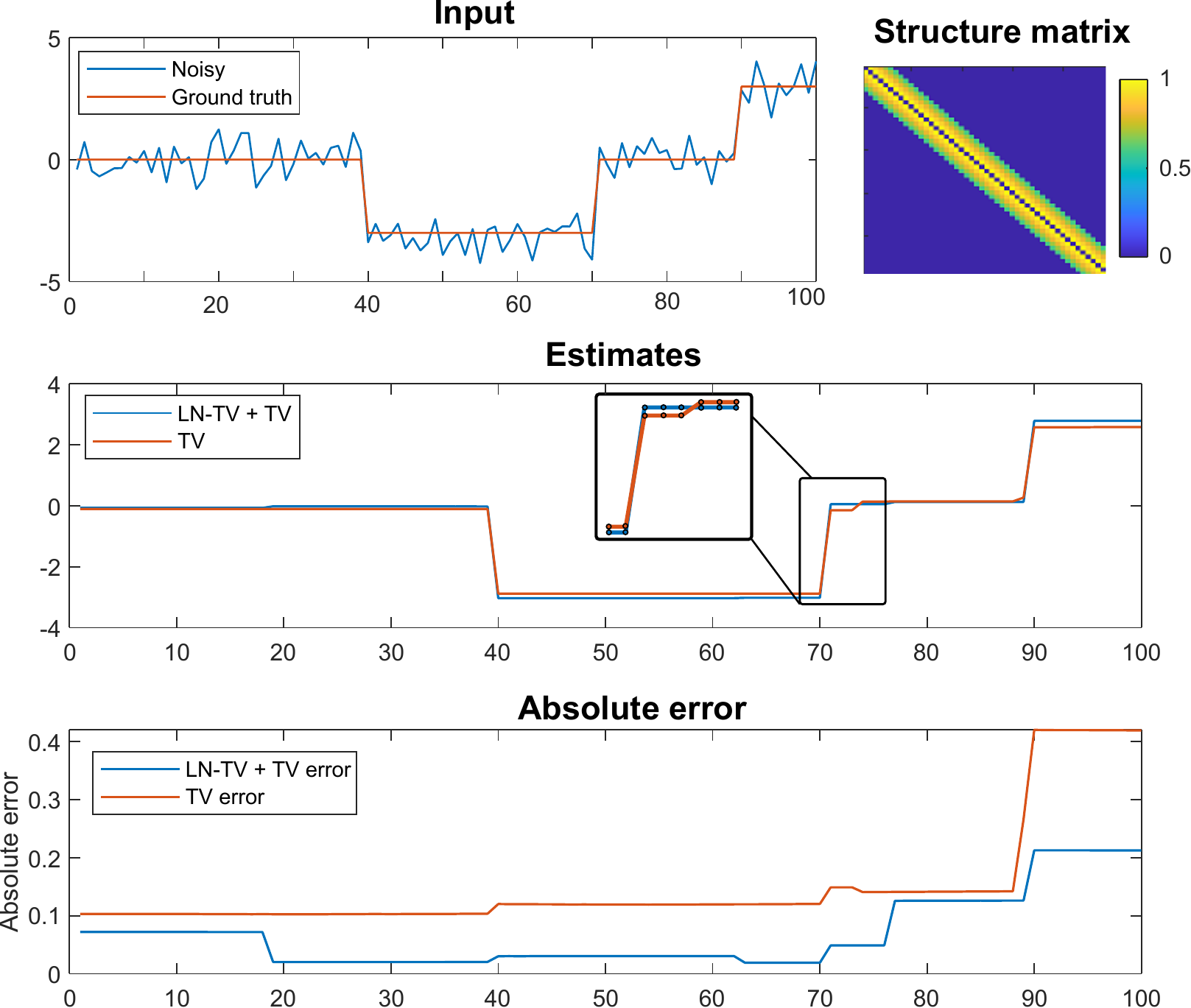}}
\caption{Combining local neighborhood total variation with standard total variation can help to reduce bias and avoid `staircasing' artefacts.}
\label{fig:LNTVcombined}
\end{figure}

For more general problems, where a piecewise approximation is helpful but the segments do not necessarily have a minimum length well known \emph{a priori}, the LN-TV penalty can be combined with the traditional TV by solving a problem of the form:
\begin{equation}
\hat{\Vecx} = \argmin_{\Vecx} \frac{1}{2}\|\VecA\Vecx - \Vecy\|_2^2 + \underbrace{\lambda_1|\VecD \Vecx| \transpose \VecS |\VecD \Vecx|}_{\rm LN-TV} + \underbrace{ \lambda_2 \|\VecD \Vecx\|_1}_{\rm TV}
\label{eqn:tvcombined}
\end{equation}
where the $\lambda$'s are parameters that control the strength of the regularization. If $\VecA = \VecI$, this acts to denoise the input $\Vecy$, but different choices of $\VecA$ can be used to solve deconvolution or deblurring problems, super-resolution, and so on.
The addition of the LN-TV term can help avoid the `staircasing' artefacts that often arise when using TV, and allow for weaker TV regularization, which reduces the bias in the result.
We present an example of these favorable effects in Fig.~\ref{fig:LNTVcombined}.
Here, the $\VecS$ matrix used promotes sparsity in $\VecD \Vecx$ within five indices either side of a nonzero component. This formulation is distinct from the one discussed in (\ref{eqn:optprob}), as instead of a constraint we are using the SWAGGER term as a regularization with a strength prescribed by $\lambda_1$. This allows us to have some more nuance by having a linear ramp from 1 to 0.5 in the band around the diagonal, to make close changes less likely than ones a greater distance away. 

We can extend this formulation to a 2D setting also. Using the penalty formulation in (\ref{eqn:tv2d}), we promote horizontal and vertical differences that occur only with some minimal separation distance:
\begin{equation}
R(\Vecx) = \lambda(\underbrace{|\VecD_h \Vecx| \transpose \VecS_h |\VecD_h \Vecx|}_{\rm horizontal} + \underbrace{|\VecD_v \Vecx| \transpose \VecS_v |\VecD_v \Vecx|}_{\rm vertical} ),
\label{eqn:tv2d}
\end{equation}
where $\VecD_v$ and $\VecD_h$ are matrices that take finite differences along the vertical columns and horizontal rows, respectively, and $\VecS_v$ and $\VecS_h$ impose local neighborhood sparsity in the vertical and horizontal directions.
We choose $\Phi(\Vecx) = |\Vecx|^\kappa$ which introduces another parameter $\kappa$ as discussed in Section~\ref{sec:phichoice}.

\begin{figure}
\centerline{\includegraphics[width=0.9\columnwidth]{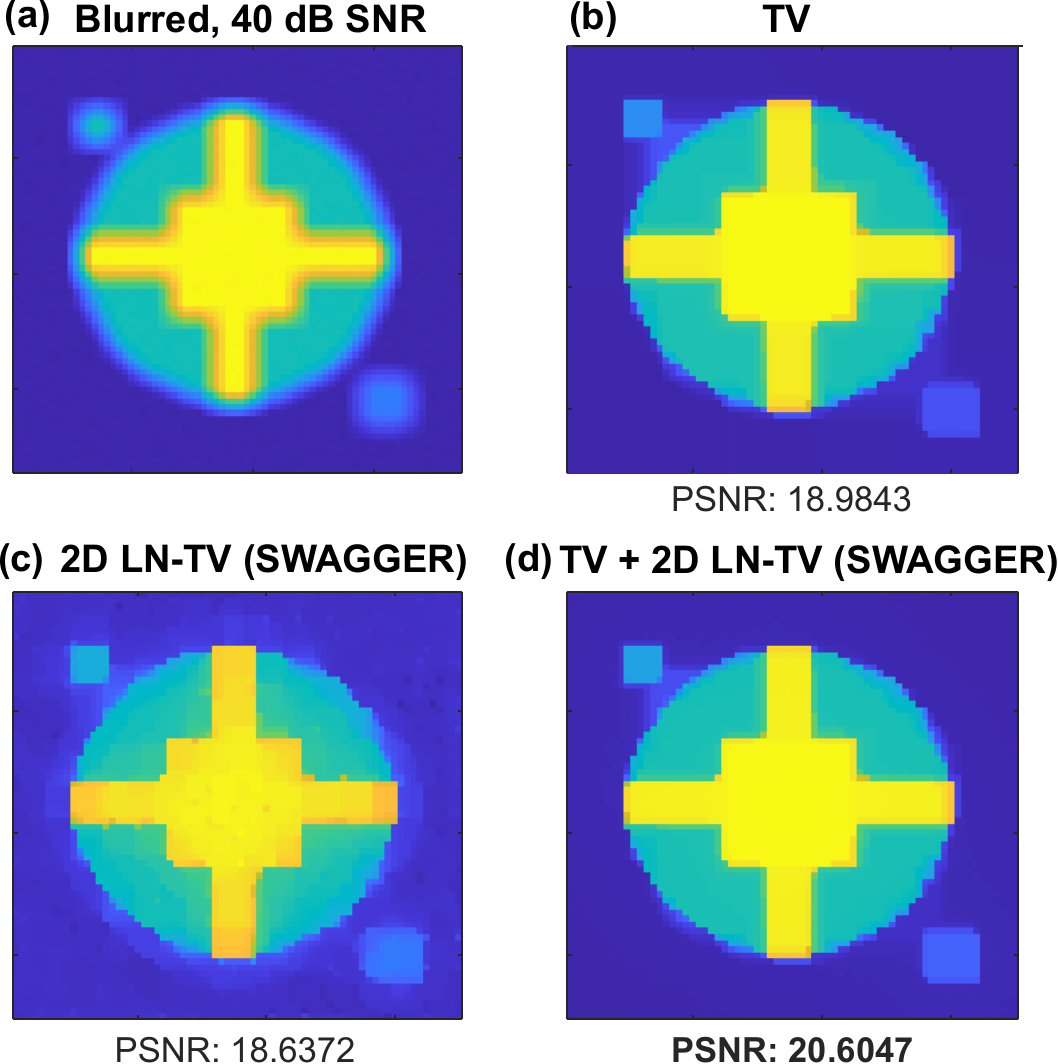}}
\caption{\textbf{(a)} A test image blurred with a Gaussian kernel and corrupted with Gaussian noise. \textbf{(b)} Reconstruction using total variation. \textbf{(c)} Reconstruction using local neighborhood TV only. \textbf{(d)} Reconstruction using a combination of both.}
\label{fig:LNTV2D}
\end{figure}

We find that, in general, $\kappa$ should be kept somewhat smaller than 1 if the initial estimate is especially blurred or noisy, and larger than 1 if there exist hard edges in the initial estimate that should be maintained.
For the results shown in Fig.~\ref{fig:LNTV2D}, $\kappa = 0.75$.

\subsection{Modeling Mutual Surface Occlusions in Non-Line-of-Sight (NLOS) Imaging}
\label{sec:erti}
The aim of NLOS imaging is to form reconstructions of scenes hidden from the direct line of sight of the observer. Typical methods to achieve this rely on the use of ultra-fast pulsed lasers and time-resolved single photon counting to infer the hidden geometry by measuring round-trip distances 
\cite{Kirmani2009, Velten2012, Heide2014, OToole2018}.
Often, to recover a 3D scene, the hidden area will be discretized into patches or surfaces in some manner.
Some of these patches may occlude others: if a surface is actually present in the scene, light coming from others behind it may have no paths to the measurement device.
As the discretization is prescribed beforehand, the surfaces occluded by each other are known \emph{a priori}, and they should be mutually exclusive.
Therefore, we can use SWAGGER to ensure that groups of surfaces that occlude one another will be one-sparse in the recovered output.
This provides two benefits. Firstly, the recovered output will be physically plausible as surfaces that cannot contribute any light in the measurements will not be present in the solution.
Secondly, ensuring this physical constraint is met \emph{during} the optimization process should lead to better estimations overall. 

This idea extends to any situation in which an ideal discretization of some physical area, field, etc., involves some mutual exclusivity between elements of the discretization.
Typically, to deal with this one may instead use a discretization that is less desirable but avoids or reduces this exclusivity requirement, use post-processing to make the recovery fit the expected structure, or simply ignore it.
Using SWAGGER can instead ensure that reconstructions both fit the expected structure prescribed by the underlying physics of the problem, and also enjoy improved estimates due to the knowledge of the structure being used within the optimization procedure.

One NLOS imaging system, Edge-Resolved Transient Imaging~\cite{ERTI}, uses a laser to illuminate the floor at numerous positions in an arc around a vertical edge, such as a doorway, to form a 3D reconstruction of the room beyond.
For each measurement position, a histogram of the arrival times of photons reaching a single photon avalanche diode (SPAD) detector is accumulated using time-correlated single photon counting.
Each subsequent laser position illuminates more of the hidden scene, and by taking differences between measurement histograms one can recover a noisy measurement containing photon arrivals originating mostly from within a single wedge in the hidden area. 

\begin{figure}
\centerline{\includegraphics[width=\columnwidth]{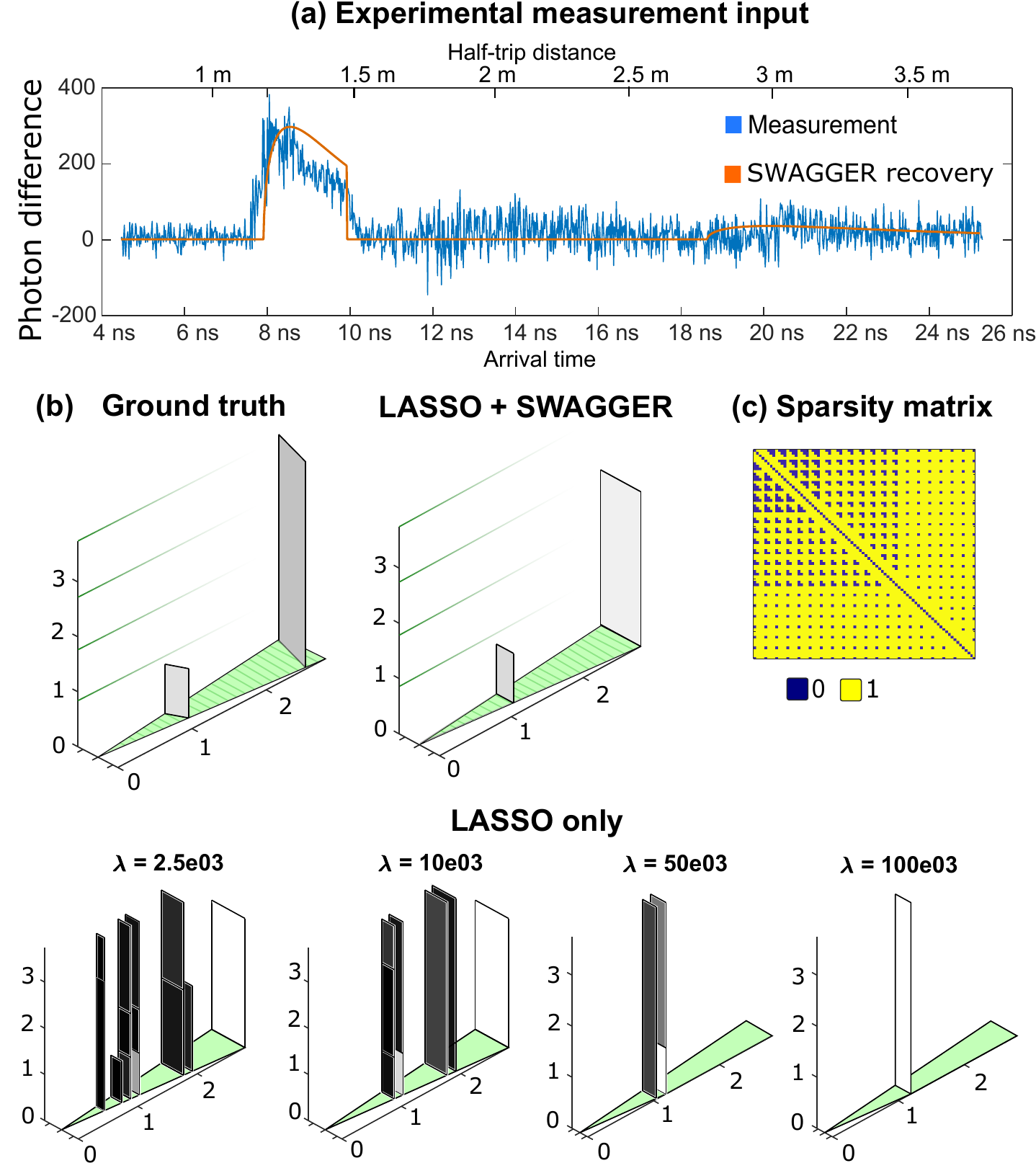}}
\caption{\textbf{(a)} Experimental measurement data containing predominantly photon returns from a single wedge in a hidden area. The data corresponds to the 18th wedge in the `Staircase' scene, Figure~4 in~\cite{ERTI}.
\textbf{(b)} The ground truth scene and the reconstruction using SWAGGER\@. The displayed reflectivity estimate is enhanced by scaling based on the amount of the surface that is occluded.
\textbf{(c)} The sparsity matrix used to model occlusion and mutual exclusivity of surfaces occupying the same space. \textbf{(d)} Reconstructions without using SWAGGER, with increasing regularization strength.
No choice of $\lambda$ gives accuracy comparable to the reconstruction using SWAGGER in (b).}
\label{fig:erti}
\end{figure}

From this difference histogram measurement, one can attempt to fit surfaces to the scene within a wedge, to reconstruct the hidden area. It is assumed that surfaces always start at the ground, extend to a certain height, and are a certain distance from the corner. We can form a reconstruction by solving the following problem:
\begin{equation}
 \begin{aligned}
\hat{\Vecx} =
& \argmin_{\Vecx \geq 0} 
& & \frac{1}{2}\|\VecA\Vecx - \Vecy\|_2^2  +  \lambda_1 \| \Vecx\|_1 \\
& \text{s.t.  }
& & |\Vecx| \transpose \VecS |\Vecx| = 0
\end{aligned}
\label{eqn:}
\end{equation}
where $\Vecy$ is a difference histogram measurement and $\VecA$ is a matrix governing the light transport to and from surfaces at different distances from the edge, and with different heights.
Fig.~\ref{fig:erti} shows reconstructions of a wedge using experimental data from~\cite{ERTI}.
One reconstruction uses only the $\|\Vecx\|_1$ term and no constraints, and a second that includes the SWAGGER constraint.
In~\cite{ERTI}, a Markov chain Monte Carlo (MCMC) method was used to form the reconstructions as it could implicitly include occlusions model and needed to only estimate parameters for a small number of surfaces, rather than use a full scene discretization.
Here, we use SWAGGER and a full scene discretization to ensure that surfaces which should occlude those behind them, and vice versa, are accurately modeled.
Secondly, it permits us to use a discretization with mutually exclusive elements---surfaces in the same position but with different heights---where only one should be present in the solution.
As such, the result is a physically plausible reconstruction, which is quite accurate in support, height, and reflectivity.

\section{Conclusion}

Using SWAGGER as a constraint allows one to recover solutions to problems that fit a known sparsity structure, including any number of one-sparse overlapping groups within the solution vector or a transform thereof. This is useful in modeling a wide variety of phenomena, such as occlusions in imaging problems. 
We illustrate the usefulness of the proposed penalty for some example use cases where, by simply constructing a sparsity structure matrix $\VecS$, one can ensure estimates fit known hard structured sparsity constraints derived from domain knowledge of the problem.
This results in solutions that are both physically plausible and also improved overall. 
Similarly, SWAGGER can be used as a regularization term with some strength, to promote solutions that fit a sparsity structure whilst allowing more nuance, which gives rise to local-neighborhood total variation. SWAGGER can be readily combined with other priors, regularization terms and constraints to improve estimations. A proximal (sub)gradient algorithm (and accelerated algorithm) have been proposed to quickly solve SWAGGER-constrained and SWAGGER-regularized problems.

\appendix

\subsection{Proof of Proposition~\ref{prop:constraints}}
\label{sec:penbounds}

Suppose both conditions on Proposition~\ref{prop:constraints} hold.
By condition~2), every entry of $\Phi(\VecB\Vecx)$ is nonnegative.
By condition~1), every entry of $\VecS$ is also nonnegative.
Therefore, each entry in the matrix triple product
$\Phi(\VecB \Vecx) \transpose \, \VecS \, \Phi(\VecB \Vecx)$
is nonnegative,
ensuring that the result is nonnegative.

\subsection{Proximity Operator Algorithm for $\ell_1^2$}
\label{sec:l12}
Algorithm~\ref{alg:prox} reproduces the method from~\cite{l12proof}
to compute
the proximal operator for the $\ell_1^2$ norm (with $O(N \log N)$ complexity).
\begin{algorithm}[H]
 \caption{Proximity operator of $\ell_1^2$}
  \label{alg:prox}
 \begin{algorithmic}[1]
 \renewcommand{\algorithmicrequire}{\textbf{Input:}}
 \renewcommand{\algorithmicensure}{\textbf{Output:}}
 \renewcommand{\algorithmicfor}{\textbf{while}}
 \REQUIRE $\Vecz \in \mathbb{R}^N$, $\lambda \geq 0$
 \ENSURE  $\hat{\Vecz} = \argmin_\Vecq \|\Vecz - \Vecq\|_2^2 + \frac{1}{2} \lambda \|\Vecq\|_1^2$
 
 \STATE sort entries of $|\Vecz|$ into $\Vecy$ s.t. ($\Vecy_1 \geq \cdots \geq \Vecy_N)$
 \STATE set $\rho = \max \{ j \in \{1...N\} \, | \, \Vecy_j - \frac{\lambda}{1+j\lambda} \sum_{r=1}^j \Vecy_r > 0 \}$

 \RETURN $\hat{\Vecz} = \max \left(|\Vecz| - \tau, 0 \right)\text{sign}(\Vecz)$ \\ where $\tau = \frac{\lambda}{1+j\lambda} \sum_{r=1}^\rho \Vecy_r$
 \end{algorithmic} 
 \end{algorithm}

\subsection{Accelerated SWAGGER Algorithm}
\label{apdx:ncapg}

In~\cite{NCAPG}, a monotonic accelerated proximal gradient method is outlined that is guaranteed to converge to a stationary point for nonconvex problems.
Algorithm~\ref{alg:swagaccel} applies this acceleration to Algorithm~\ref{alg:swag},
where $\alpha_x$, $\alpha_y$ and $\alpha_\lambda$ are step sizes that can be fixed or computed dynamically with a back-tracking scheme. The convergence proof requires that $f(\Vecx)$ is a proper function with Lipschitz continuous gradients (e.g., the least squares loss), and the SWAGGER penalty is proper and lower semicontinuous, which is the case for $\Phi(\Vecx) = |\Vecx|$. 

 \begin{algorithm}[H]
 \caption{Accelerated proximal subgradient method for SWAGGER}
  \label{alg:swagaccel}
 \begin{algorithmic}[1]
 \renewcommand{\algorithmicrequire}{\textbf{Input:}}
 \renewcommand{\algorithmicensure}{\textbf{Output:}}
 \renewcommand{\algorithmicfor}{\textbf{while}}
 \REQUIRE $\Vecz_1 = \Vecx_1 = \Vecx_0$, $t_1 = 1$, $t_0 = 0$, $\overline{\VecS} = \mathbbm{1}\mathbbm{1}\transpose - \VecS $.
 \ENSURE  $\hat{\Vecx}$

  \FOR {not converged}
  \STATE $\Vecy^k = \Vecx^k + \frac{t^{k-1}}{t^k}(\Vecz^k - \Vecx^k) + \frac{t^{k-1}-1}{t^k}(\Vecx^k - \Vecx^{k-1})$
  \STATE $\Vecz^{k+1} = P_{\alpha_y \lambda^k}( \Vecy^{k} - \alpha_y(\nabla_y ( f(\Vecy^k) - \lambda^k|\Vecy^k|\transpose \, \overline{\VecS} \, |\Vecy^k|))$
  \STATE $\Vecv^{k+1} = P_{\alpha_x \lambda^k}( \Vecx^{k} - \alpha_x(\nabla_x (f(\Vecx^k) - \lambda^k|\Vecx^k|\transpose \, \overline{\VecS} \, |\Vecx^k|))$

\STATE  $\Vecx^{k+1} = \begin{cases} 
          \Vecz^{k+1}, & F(\Vecz^{k+1}) \leq F(\Vecv^{k+1}); \\
          \Vecv^{k+1}, & \text{otherwise}
       \end{cases}$

  \STATE $t^{k+1} = \frac{1}{2}(\sqrt{4(t^k)^2 + 1} + 1)$
  
  \STATE $\lambda^{k+1} = \lambda^k + \alpha_\lambda(|\Vecx^{k}| \transpose \VecS |\Vecx^{k}|)$
  \ENDFOR
 \RETURN $\hat{\Vecx} = \Vecx^{k+1} $ 
 \end{algorithmic} 
 \end{algorithm}

\subsection{Generating Test Vectors}
\label{sec:genx}
For the experiments in Section~\ref{sec:compare}, we must generate a test vector $\Vecx$ for each trial which fits the sparsity structure described by a specific $\VecS$ matrix.
To do so we use Algorithm~\ref{alg:genx},
where in Step 4,
$\Vecu_j$ is drawn from the continuous uniform distribution on $[0.5,1.5]$.

 \begin{algorithm}[H]
 \caption{Generate a vector $\Vecx_S$ s.t. $|\Vecx^i|\transpose \VecS |\Vecx^i| \approx 0$}
  \label{alg:genx}
 \begin{algorithmic}[1]
 \renewcommand{\algorithmicrequire}{\textbf{Input:}}
 \renewcommand{\algorithmicensure}{\textbf{Output:}}
 \renewcommand{\algorithmicfor}{\textbf{while}}
 \REQUIRE $\Vecx^0 \sim N(0, 1)$, $\VecS$, $\mu \ll 1$
 \ENSURE  $\Vecx_S$

  \FOR { $|\Vecx^i|\transpose \VecS |\Vecx^i| \geq \mu$}
  \STATE $\Vecx^{i+1} = \argmin_\Vecx \|\Vecx^{i} - \Vecx\| + |\Vecx|\transpose \VecS |\Vecx|$

  \ENDFOR
 \RETURN $\Vecx_{S,j} = \begin{cases} 
          0, & \Vecx^{i+1}_j = 0 ; \\
          \text{sign}(\Vecx^{i+1}_j) \Vecu_j, & \Vecx^{i+1}_j \neq 0
       \end{cases}$
 \end{algorithmic}
 \end{algorithm}

 \subsection{The Choice of $c$ for Convex Nonconvexity}
 \label{sec:cnc}
 
The Hessian of the cost function
\begin{equation*}
  \hat{\Vecx} = \argmin_{\Vecx} \frac{1}{2}\|\VecA\Vecx - \Vecy\|_2^2 + \lambda|\Vecx| \transpose \VecS |\Vecx|
\end{equation*}
is given by $\VecA \transpose \VecA + 2\lambda \VecX\VecS\VecX$
(where $\VecX = \text{diag}(\text{sgn}(\Vecx))$). We introduce a term $c \VecI$ to the $\VecS$ matrix:
\begin{eqnarray*}
  \VecA \transpose \VecA + 2\lambda \VecX(\VecS - c \VecI)\VecX
  & = & \VecA \transpose \VecA + 2\lambda \VecX\VecS \VecX - 2\lambda c \VecX\VecI\VecX \\
  & = & \VecA \transpose \VecA + 2\lambda \VecX\VecS \VecX  - 2\lambda c \VecI.
\end{eqnarray*}
To maintain convexity, we wish to have
\begin{equation*}
   \lambda_{\text{min}}(\VecA \transpose \VecA + 2\lambda \VecX\VecS \VecX) - 2\lambda c \geq 0,
\end{equation*}
or equivalently
\begin{equation*}
   2\lambda c \leq \lambda_{\text{min}}(\VecA \transpose \VecA + 2\lambda \VecX\VecS \VecX).
\end{equation*}

Since $\VecA \transpose \VecA$ and $2\lambda \VecX\VecS \VecX$ are both
Hermitian,
\begin{equation*}
\lambda_{\text{min}}(\VecA \transpose \VecA + 2\lambda \VecX\VecS \VecX)  \geq \lambda_{\text{min}}(\VecA \transpose \VecA) + \lambda_{\text{min}}(2\lambda \VecX\VecS \VecX)
\end{equation*}
by Weyl's inequality~\cite{weyl}.
Finally, we note
\begin{equation*}
   \min_\VecX \lambda_{\text{min}}(2\lambda \VecX\VecS \VecX)) = 2\lambda \, \lambda_\text{min}(\VecS),
\end{equation*}
providing us with a bound for the choice of $c$:
\begin{equation*}
c \leq \frac{1}{2\lambda}\lambda_{\text{min}}(\VecA \transpose \VecA) + \lambda_{\text{min}}(\VecS ).
\end{equation*}

\section*{Acknowledgment}

We thank Y. Ma for comments on the manuscript, and S. W. Seidel, J. Rapp, J. Murray-Bruce and I. Selesnick for discussions.

\bibliography{lit}

\begin{thebibliography}{10}
\providecommand{\url}[1]{#1}
\csname url@samestyle\endcsname
\providecommand{\newblock}{\relax}
\providecommand{\bibinfo}[2]{#2}
\providecommand{\BIBentrySTDinterwordspacing}{\spaceskip=0pt\relax}
\providecommand{\BIBentryALTinterwordstretchfactor}{4}
\providecommand{\BIBentryALTinterwordspacing}{\spaceskip=\fontdimen2\font plus
\BIBentryALTinterwordstretchfactor\fontdimen3\font minus
  \fontdimen4\font\relax}
\providecommand{\BIBforeignlanguage}[2]{{%
\expandafter\ifx\csname l@#1\endcsname\relax
\typeout{** WARNING: IEEEtran.bst: No hyphenation pattern has been}%
\typeout{** loaded for the language `#1'. Using the pattern for}%
\typeout{** the default language instead.}%
\else
\language=\csname l@#1\endcsname
\fi
#2}}
\providecommand{\BIBdecl}{\relax}
\BIBdecl

\bibitem{LASSO}
R.~Tibshirani, ``Regression shrinkage and selection via the lasso,'' \emph{J.
  Royal Statistical Society. Series B (Methodological)}, vol.~58, no.~1, pp.
  267--288, 1996.

\bibitem{GLASSOthesis}
S.~Bakin, ``Adaptive regression and model selection in data mining problems,''
  Ph.D. dissertation, Australian National University, 1999.

\bibitem{GLASSO}
M.~Yuan and Y.~Lin, ``Model selection and estimation in regression with grouped
  variables,'' \emph{J. Royal Statistical Society Series B}, vol.~68, no.~1,
  pp. 49--67, Feb. 2006.

\bibitem{NNsparsity}
J.~Yoon and S.~J. Hwang, ``Combined group and exclusive sparsity for deep
  neural networks,'' in \emph{Proc. 34th Int. Conf. Machine Learning}, ser.
  Proceedings of Machine Learning Research, D.~Precup and Y.~W. Teh, Eds.,
  vol.~70, 06--11 Aug 2017, pp. 3958--3966.

\bibitem{computervision1}
F.~Song, X.~Tan, and S.~Chen, ``Exploiting relationship between attributes for
  improved face verification,'' \emph{Computer Vision and Image Understanding},
  vol. 122, pp. 143--154, 2014.

\bibitem{computervision2}
A.~Wang, J.~Cai, J.~Lu, and T.-J. Cham, ``Modality and component aware feature
  fusion for {RGB-D} scene classification,'' in \emph{Proc. IEEE Conf. Computer
  Vision and Pattern Recognition (CVPR)}, Jun. 2016.

\bibitem{medicineEEG}
D.~Paz-Linares, M.~Vega-Hernández, P.~A. Rojas-López, P.~A.
  Valdés-Hernández, E.~Martínez-Montes, and P.~A. Valdés-Sosa, ``Spatio
  temporal {EEG} source imaging with the hierarchical {B}ayesian elastic net
  and elitist lasso models,'' \emph{Frontiers in Neuroscience}, vol.~11, p.
  635, 2017.

\bibitem{medicinefMRI}
N.~Rao, C.~Cox, R.~Nowak, and T.~T. Rogers, ``Sparse overlapping sets lasso for
  multitask learning and its application to {fMRI} analysis,'' in
  \emph{Advances in Neural Information Processing Systems 26}, C.~J.~C. Burges,
  L.~Bottou, M.~Welling, Z.~Ghahramani, and K.~Q. Weinberger, Eds., 2013, pp.
  2202--2210.

\bibitem{ELASSO1}
M.~Kowalski, ``Sparse regression using mixed norms,'' \emph{Applied and
  Computational Harmonic Analysis}, vol.~27, no.~3, pp. 303--324, 2009.

\bibitem{ELASSO2}
M.~Kowalski and B.~Torrésani, ``Sparsity and persistence: Mixed norms provide
  simple signal models with dependent coefficients,'' \emph{Signal Image and
  Video Processing}, vol.~3, pp. 251--264, Sep. 2009.

\bibitem{ELASSO3}
Y.~Zhou, R.~Jin, and S.~C. Hoi, ``Exclusive lasso for multi-task feature
  selection,'' in \emph{Proc. Thirteenth Int. Conf. Artificial Intelligence and
  Statistics}, ser. Proceedings of Machine Learning Research, Y.~W. Teh and
  M.~Titterington, Eds., vol.~9, 13--15 May 2010, pp. 988--995.

\bibitem{SWAG1}
\.{I}lker {Bayram} and S.~Bulek, ``A penalty function promoting sparsity within
  and across groups,'' \emph{IEEE Trans. Signal Process.}, vol.~65, no.~16, pp.
  4238--4251, Aug 2017.

\bibitem{SWAG2}
\.{I}lker {Bayram}, ``Sparsity within and across overlapping groups,''
  \emph{IEEE Signal Processing Letters}, vol.~25, no.~2, pp. 288--292, Feb
  2018.

\bibitem{eigbound}
S.~Walker and P.~{Van Mieghem}, ``On lower bounds for the largest eigenvalue of
  a symmetric matrix,'' \emph{Linear Algebra and its Applications}, vol. 429,
  no.~2, pp. 519 -- 526, 2008.

\bibitem{hollowmatrix}
Z.~B. Charles, M.~Farber, C.~R. Johnson, and L.~Kennedy-Shaffer, ``Nonpositive
  eigenvalues of hollow, symmetric, nonnegative matrices,'' \emph{SIAM J.
  Matrix Analysis and Applications}, vol.~34, no.~3, pp. 1384--1400, 2013.

\bibitem{Selesnick2017}
I.~Selesnick, ``Sparse regularization via convex analysis,'' \emph{IEEE Trans.
  Signal Process.}, vol.~65, no.~17, pp. 4481--4494, 2017.

\bibitem{l12proof}
A.~F.~T. Martins, N.~Smith, E.~Xing, P.~Aguiar, and M.~Figueiredo, ``Online
  learning of structured predictors with multiple kernels,'' in \emph{Proc.
  Fourteenth Int. Conf. Artificial Intelligence and Statistics}, G.~Gordon,
  D.~Dunson, and M.~Dudík, Eds., vol.~15, 11--13 Apr 2011, pp. 507--515.

\bibitem{NCAPG}
H.~Li and Z.~Lin, ``Accelerated proximal gradient methods for nonconvex
  programming,'' in \emph{Advances in Neural Information Processing Systems
  28}, C.~Cortes, N.~D. Lawrence, D.~D. Lee, M.~Sugiyama, and R.~Garnett,
  Eds.\hskip 1em plus 0.5em minus 0.4em\relax Curran Associates, Inc., 2015,
  pp. 379--387.

\bibitem{admm}
\BIBentryALTinterwordspacing
S.~Boyd, N.~Parikh, E.~Chu, B.~Peleato, and J.~Eckstein, ``Distributed
  optimization and statistical learning via the alternating direction method of
  multipliers,'' \emph{Found. Trends Mach. Learn.}, vol.~3, no.~1, p. 1–122,
  Jan. 2011. [Online]. Available: \url{https://doi.org/10.1561/2200000016}
\BIBentrySTDinterwordspacing

\bibitem{pshrinkage}
J.~Woodworth and R.~Chartrand, ``Compressed sensing recovery via nonconvex
  shrinkage penalties,'' \emph{Inverse Problems}, vol.~32, no.~7, p. 075004,
  May 2016.

\bibitem{MTV}
I.~{Selesnick}, ``Total variation denoising via the {M}oreau envelope,''
  \emph{IEEE Signal Process. Lett.}, vol.~24, no.~2, pp. 216--220, 2017.

\bibitem{Kirmani2009}
A.~Kirmani, T.~Hutchison, J.~Davis, and R.~Raskar, ``Looking around the corner
  using transient imaging,'' in \emph{Proc. IEEE Int. Conf. Computer Vision},
  2009, pp. 159--166.

\bibitem{Velten2012}
A.~Velten, T.~Willwacher, O.~Gupta, A.~Veeraraghavan, M.~G. Bawendi, and
  R.~Raskar, ``Recovering three-dimensional shape around a corner using
  ultrafast time-of-flight imaging,'' \emph{Nat. Commun.}, vol.~3, 2012.

\bibitem{Heide2014}
F.~Heide, L.~Xiao, W.~Heidrich, and M.~B. Hullin, ``Diffuse mirrors: {3D}
  reconstruction from diffuse indirect illumination using inexpensive
  time-of-flight sensors,'' in \emph{Proc. IEEE Conf. Computer Vision and
  Pattern Recognition}, 2014, pp. 3222--3229.

\bibitem{OToole2018}
M.~O'Toole, D.~B. Lindell, and G.~Wetzstein, ``Confocal non-line-of-sight
  imaging based on the light-cone transform,'' \emph{Nature}, vol. 555, pp.
  338--341, Mar. 2018.

\bibitem{ERTI}
J.~Rapp, C.~Saunders, J.~Tachella, J.~Murray-Bruce, Y.~Altmann, J.-Y.
  Tourneret, S.~McLaughlin, R.~M.~A. Dawson, F.~N.~C. Wong, and V.~K. Goyal,
  ``Seeing around corners with edge-resolved transient imaging,''
  arXiv:2002.07118v1 [eess.IV]., Feb. 2020.

\bibitem{weyl}
R.~A. Horn and C.~R. Johnson, \emph{Matrix Analysis: Second Edition}.\hskip 1em
  plus 0.5em minus 0.4em\relax Cambridge University Press, 2012.

\end{thebibliography}
\bibliographystyle{IEEEtran}

\end{document}